\documentclass[12pt,letterpaper,twoside]{amsart}

\usepackage[ngerman]{babel}
\usepackage{graphicx}
\usepackage{latexsym}
\usepackage{amsmath,amssymb}
\usepackage{ifpdf,stmaryrd}

\usepackage[DIV14,BCOR2cm]{typearea}
\theoremstyle{definition}
\newtheorem{Definition}{Definition}
\newtheorem{Theorem}[Definition]{Theorem}
\newtheorem{Corollary}[Definition]{Corollary}
\newtheorem{Lemma}[Definition]{Lemma}
\newenvironment{Beweis}{\noindent{Proof:}}{%
  \hspace*{\fill}$\Box$\par\vskip2ex}

\begin{document}


\title[Applications of the Graph Minor Theorem to algebraic structures I]{Applications of the Graph Minor Theorem to algebraic structures I}
\author{Tobias Ahsendorf}
\curraddr{Fakult\"at f\"ur Mathematik, Universit\"at Bielefeld, P.O.Box 100 131, D-33501 Bielefeld, Germany}
\email{tobias.ahsendorf@gmail.com}
\maketitle

\section*{Abstract}
We use a correspondence between Graphs\footnote{Graphs are always finite in this paper.} and finite subsets of certain semigroups (the factorial and commutative ones, for example $R \backslash \{0\}$ if $R$ is a unique factorization domain) and the Graph Minor Theorem to get a characterization of infinite sequences of finite subsets of those semigroups. Especially interesting is the case if $\mathbb{N}$\footnote{Here $\mathbb{N}$ starts at $1$.} is the regarded semigroup. 

\section{Some aspects of Graph (Minor) Theory}
Robertson and Seymour proved the so-called Graph Minor Theorem in twenty groundbreaking papers, published from 1983 till 2004. To state this theorem, we need to introduce the definition of a minor.

\begin{Definition}(see \cite{RDE}) \label{min}
Let $G$ and $H$ be graphs. We say $H$ is a \emph{minor} of $G$, in notation $H \preccurlyeq G$, when there is a partition $\{V_h\}_{h \in V(H)}$ of a subset of the vertex set of $G$ with $G[V_h]$ connected for each $h$, such that there is a $V_h$-$V_{h'}$ edge if $hh' \in E(H)$. 
\end{Definition}
This definition doesn't stand there directly, for this you have to read a bit between the lines.
An equivalent definition, which can maybe better be used to understand the term of a minor is the following, but it won't be necessary for the rest of the paper.
\begin{Definition}(see \cite{LL}) $H \preccurlyeq G$ if $H$ can be obtained from $G$ by a sequence of the following operations:
\begin{itemize}
\item Delete an edge
\item Delete a vertex
\item Contract an edge
\end{itemize}
The \emph{contraction} of an edge is defined in the obvious way.
\end{Definition}

It is easy to see that the minor relation is reflexive and transitive in the class of finite graphs. Now we can state the Graph Minor Theorem.
\begin{Theorem}(see \cite{RSXX}) \label{RS}The finite graphs are well-quasi-ordered by the minor relation $\preccurlyeq$, i.e. $\preccurlyeq$ is reflexive and transitive and for an infinite sequence of graphs $G_1,G_2,\dots$ there are $i,j$ with $i < j$ and $G_i \preccurlyeq G_j$.
\end{Theorem}

We actually can show an even stronger argument. For this we have to use a Ramsey-like Theorem, in fact we can deduce the classical Ramsey Theorem from it, so it can be regarded as a more general version of it. Therefore let $X$ be a set. A $c$\emph{-colouring} is a partition into $c$ classes. With $[X]^k$ we denote the set of all subsets with $k$ elements. We say that $Y \subseteq X$ is \emph{monochromatic} under a certain $c$-colouring of $[X]^k$, if $[Y]^k$ has the same colour, which means that all elements belong to the same partion class.

\begin{Theorem}(see \cite{RDE}) \label{Ramsey} Let $k,c$ be positive integers and $X$ an infinite subset. If $[X]^k$ is coloured with $c$ colours, then there exists a monochromatic $Y \subseteq X$ with infinitely many elements.
\end{Theorem}

With the help of this theorem, we are able to show that in any infinite sequence of graphs $G_1,G_2,\dots$ there is not only one pair $(i,j)$ with $i<j$, such that $G_i \preccurlyeq G_j$, but there is an infinite sequence $\{i_j  \}_{j \in \mathbb{N}}$ of positve integers such that $i_j < i_{j+1}$ and $G_{i_j} \preccurlyeq G_{i_{j+1}}$ for all $j$ and so because of the fact that the minor relation is transitive $G_{i_j} \preccurlyeq G_{i_k}$ if $j \leq k$. If these conditions hold we call $G_{i_1},G_{i_2},\dots$ an \emph{increasing} subsequence of $G_1,G_2,\dots$ under $\preccurlyeq$.

\begin{Theorem}(see \cite{RDE}) \label{inc} Every infinite sequence of graphs contains an infinite in\-creasing subsequence under $\preccurlyeq$.
\end{Theorem}

\begin{Beweis}
Let $G_1,G_2,\dots$ be a infinite sequence of graphs. Now we take $[\mathbb{N}]^2$ and for all elements $\{i,j\}$ of this set, without loss of generality $i<j$, we colour $\{i,j\}$ with red if $G_j \preccurlyeq G_i$ and $G_j \neq G_i$, with yellow if $G_i$ and $G_j$ are incomparable and with green if $G_i \preccurlyeq G_j$. Now we have a $3$-colouring of an infinite set and we can use Theorem \ref{Ramsey}, so we obtain a infinite subset $Y$ of $\mathbb{N}$ which is monochromatic. Now suppose that this colour is red or yellow. This would contradict Theorem \ref{RS}, because $\{G_y \}_{y \in Y}$ would be a counterexample. The only colour that remains is green, so $\{G_y \}_{y \in Y}$ is an infinte increasing subsequencce under $\preccurlyeq$.
\end{Beweis}

\section{Factorial commutative semigroups}
First of all we have to define what we mean with the term of a factorial semigroup.
\begin{Definition}
Let $S =(S,e,\cdot)$ be a commutative semigroup (here a semigroup has a neutral element $e$). We say $S$ is \emph{factorial} if for each $s \in S$ we can write it in the unique equation $s = \varepsilon \cdot p_1 \cdots p_r$ with $r \in \mathbb{N}\, \cup \{0\}$, $\varepsilon \in S^{\times} := \{s \in S | \exists s'\in S: s \cdot s' = e \}$ and prime elements $p_1,\dots,p_r \in P(S) := \{ p \in S \backslash S^{\times} | \forall a,b \in S: p|a\cdot b \Rightarrow p|a \textnormal{ or }p|b \}$. With $c|d$ and $c,d \in S$ we mean, that there is $f \in S$ with $fc=d$. We denote the class of all such semigroups with $FCSG$ (factorial commutative semigroup).
\end{Definition}

The most important of such semigroups is $\mathbb{N} = (\mathbb{N},1,\cdot)$. 

\begin{Lemma} \label{cs}
If $S =(S,e,\cdot)$ is a commutative semigroup (where $S^{\times}$ is defined in the similar way for nonfactorial ones as in the previous definition), then $a,b \in S^{\times} \Leftrightarrow ab \in S^{\times}$ or equivalently:  $a$ or $b \in S \backslash S^{\times} \Leftrightarrow ab \in S\backslash S^{\times}$. \\
This means (among others), if $S \in FCSG$, then any non-empty product of prime elements is an element of $S \backslash S^{\times}$.
\end{Lemma}
\newpage
\begin{Beweis} First of all we prove the first statement.
\begin{itemize}
\item["'$\Rightarrow$"'] Let $c,d \in S$ with $ac=e$ and $bd=e$. Then $(ab)(cd)=(ac)(bd) = e \cdot e = e$, so $ab \in S^{\times}$. 
\item["'$\Leftarrow$"'] $ab \in S^{\times}$, so there exists $c \in S$ with $e =abc=a(bc)=b(ac)$, i.e. $a,b \in S^{\times}$.
\end{itemize}
Now we show by induction on the number $n$ of prime elements occuring in a certain product, that the last statement is true. \\
If $n=1$ our assumption is clear by definition. Now we assume that this is true for a certain $n$ and we show that this so for $n+1$ as well. \\
Let $a = \prod_{i=1}^{n} p_i$ be a product of $n$ prime elements, so $a \in S\backslash S^{\times}$. If $b=p_{n+1}$ is any element of $P(S) \subseteq S\backslash S^{\times}$, so we know by the first part of this Lemma that $ab= \prod_{i=1}^{n+1} p_i \in S \backslash S^{\times}$. 
\end{Beweis}
The next definition generalizes the idea of the greatest common divisor from $\mathbb{N}$ to all elements of $FCSG$.

\begin{Definition} \label{gcd}
Let $S=(S,e,\cdot) \in FCSG$ and $a_1,a_2 \in S$ be arbitrary elements with $a_i = \varepsilon_i \cdot \prod_{j=1}^{t} p_j^{\alpha_i(j)}$ for $i=1,2$ and with $t, \alpha_1(1),\dots,\alpha_1(t),\alpha_2(1),\dots,\alpha_2(t) \in \mathbb{N}\cup \{0\}$, $p_1,\dots,p_t \in P(S)$ and $p_i \neq p_j$ if $i \neq j$. Such a choose is possible - maybe with many $\alpha_i(j) =0$. The greatest common divisor of $S$, denoted by $\gcd_S$, is defined as $\gcd_S : S \times S \rightarrow S$ with $\gcd_S(a_1,a_2) := \prod_{j=1}^{t} p_j^{\min \{\alpha_1(j),\alpha_2(j) \}}$, where $\prod_{j=1}^{t} p_j^{0}:= e$ for all $t \in \mathbb{N} \, \cup \{0\}$ and $p_1,\dots,p_t \in P(S)$.
\end{Definition}
It is clear, that $\gcd_S$ is welldefined for all $S \in FCSG$ because of the unique factorization of every element of $S$.
Now we prove a Lemma which shows a certain invariance of the $\gcd_S$.
\begin{Lemma} \label{inv} Let $a_1,a_2,a_3,a_4 \in S \in FCSG$. If $\gcd_S(a_1,a_2) \notin S^{\times}$, then \\ $\gcd_S(a_1a_3,a_2a_4) \notin S^{\times}$. 
\end{Lemma}

\begin{Beweis}
Let $a_i$ be in the same representation like in Definition \ref{gcd}, i.e. $a_i = \varepsilon_i \cdot \prod_{j=1}^{t} p_j^{\alpha_i(j)}$ for $i=1,2,3,4$. Now $\gcd_S(a_1,a_2) = \prod_{j=1}^{t} p_j^{\min \{\alpha_1(j),\alpha_2(j) \}}$ and $\gcd_S(a_1a_3,a_2a_4) = \prod_{j=1}^{t} p_j^{\min \{ \alpha_1(j)+\alpha_3(j), \alpha_2(j)+\alpha_4(j) \}}$, we can deduce, that $\min \{ \alpha_1(j)+\alpha_3(j), \alpha_2(j)+\alpha_4(j) \} \geq \min \{\alpha_1(j),\alpha_2(j) \}$. Because we know that there is $j$ with $\min \{\alpha_1(j),\alpha_2(j) \} >0$ (otherwise $\gcd_S(a_1,a_2)=e \in S^{\times}$), we know that $\gcd_S(a_1a_3,a_2a_4)$ is a non-empty product of prime elements and by Lemma \ref{cs} we have proved this Lemma.
\end{Beweis}

\section{The combination of Graph Minor Theory and FCSGs}

With the techniques as they are brought together we are now able to formulate the Theorem which there was refered to nebulously in the abstract. It actually consists of two assertions of which one is stronger than the other one, but the weaker one is easier to formulate and to understand. For $M \subseteq S=(S,e,\cdot) \in \mathbb{N}$ we define $m(M) := \prod_{x \in M} x$ if $M \neq \varnothing$ and $m(\varnothing) := e$.

\newpage
\begin{Theorem}  \label{the}
Let $S=(S,e,\cdot)\in FCSG$ and $\{M_i\}_{i \in \mathbb{N}}$ be an infinte sequence of finite subsets of $S$. Then there exists a sequence $\{i_j\}_{j \in \mathbb{N}}$ of positive integers with $i_j < i_{j+1}$ for all $j$ with the following properties (for each $j$):
\begin{itemize}
\item We can partitionize a subset of $M_{i_{j+1}}$ in $\{M_{i_{j+1},k}\}_{k \in M_{i_j}}$ with the following conditions for each $k,k' \in M_{i_j}$ with $k\neq k'$: 
\begin{itemize}
\item If $\gcd_S(k,k') \notin S^{\times}$, then $\gcd_S(m(M_{i_{j+1},k}),m(M_{i_{j+1},k'})) \notin S^{\times}.$ 
\item For each $x,y \in M_{i_{j+1},k}$ there is $l \in \mathbb{N} \, \cup \{0\}$ and $a_1,\dots,a_l \in M_{i_{j+1},k}$ such that $\gcd_S(x,a_1),\gcd_S(a_1,a_2),\dots,\gcd_S(a_{l-1},a_l),\gcd_S(a_l,y) \notin S^{\times}$.
\end{itemize}
\item We can partitionize  $M_{i_{j+1}}$ in $\{M'_{i_{j+1},k}\}_{k \in M_{i_j}}$ with the above conditions, only the second condition does not necessarily hold for one (arbitrary) $k_0 \in M_{i_j}$.
\end{itemize} 
\end{Theorem}

\begin{Beweis} First we show the first assertion, the second will follow immediatly, as we will see.
For each $i$ we define $G_i :=(M_i, \{ \{x,y\}\in [M_i]^2 | \gcd_S(x,y) \notin S^{\times} \})$, because $M_i$ is the vertex set it es clear that $G_i$ is a (finite) graph. Now regard the sequence $G_0,G_1,\dots$. By Theorem \ref{inc} we get a sequence $\{i_j\}_{j \in \mathbb{N}}$ of positive ascending integers with $$H:=G_{i_j} \preccurlyeq G_{i_{j+1}}=:G.$$ Let $\{V_h\}_{h \in V(H)}$ be like in Definition \ref{min}. Because of $M_{i_j} = V(H)$ and $M_{i_{j+1}}=V(G)$ we can set $M_{i_{j+1},k}:= V_k$. Because the $V_k$ are connected in $G$ and by definition when there is an edge between two elements of $M_i$ in $G_i$ it follows the first part of the first assertion. Hence we only have to show, that if $\gcd_S(k,k')\notin S^{\times}$ and $k \neq k'$, then $\gcd_S(m(M_{i_{j+1},k}),m(M_{i_{j+1},k'})) \notin S^{\times}$. When $\gcd_S(k,k') \notin S^{\times}$ then $kk' \in E(H)$, there is by definition a $V_k$-$V_{k'}$ edge, i.e. there is $x \in M_{i_{j+1},k}$ and $x' \in M_{i_{j+1},k'}$ such that $\gcd_S(x,x')\notin S^{\times}$, but then by Lemma \ref{inv} we get $\gcd_S(x \cdot m(M_{i_{j+1},k} \backslash \{x\}), x' \cdot m(M_{i_{j+1},k'} \backslash \{x'\}))= \gcd_S(m(M_{i_{j+1},k}),m(M_{i_{j+1},k})) \notin S^{\times}$, and so the entire first part of this Theorem is proved. \\
Now have a look at the second assertion (we simply overtake the notation and main idea of the first part of this proof): If $M_{i_j}=\varnothing$ then the fact is trivial. So we may assume that $M_{i_j} \neq \varnothing$. Let $M$ be defined as  $M_{i_{j+1}} \backslash ( \bigcup_{k \in M_{i_j}} M_{i_{j+1},k} )$. Choose an arbitraty, but fixed $k_0 \in M_{i_j}$ and define $M'_{i_{j+1},k}:=M_{i_{j+1},k}$ if $k \neq k_0$ and $M'_{i_{j+1},k_0}:=M_{i_{j+1},k_0} \cup M$. It is clear, that this is a partition of $M_{i_{j+1}}$. So there is only nontrivial fact we have to show: For each $k,k' \in M_{i_j}$ with $\gcd_S(k,k') \notin S^{\times}$ and $k \neq k'$, we have $\gcd_S(m(M'_{i_{j+1},k}),m(M'_{i_{j+1},k'})) \notin S^{\times}.$ If $k,k' \neq k_0$ nothing has changed, so this consequence is true. If one of these elements equals $k_0$, say $k$, then we get the result by Lemma \ref{inv}, since $\gcd_S(k_0,k')\notin S^{\times} \Rightarrow \gcd_S(m(M_{i_{j+1},k_0}),m(M_{i_{j+1},k'})) \notin S^{\times} \Rightarrow \gcd_S(m(M'_{i_{j+1},k_0}),m(M'_{i_{j+1},k'})) \notin S^{\times}$, where we have used Lemma \ref{inv} in the last implication.
\end{Beweis}

As a Corollary we can spezialize this Theorem on $\mathbb{N}$.
\begin{Corollary}
Let $\{M_i\}_{i \in \mathbb{N}}$ be an infinte sequence of finite subsets of $\mathbb{N}$. Then there exists a sequence $\{i_j\}_{j \in \mathbb{N}}$ of positive integers with $i_j < i_{j+1}$ for all $j$ with the following properties (for each $j$):
\begin{itemize}
\item We can partitionize a subset of $M_{i_{j+1}}$ in $\{M_{i_{j+1},k}\}_{k \in M_{i_j}}$ with the following condition for each $k,k' \in M_{i_j}$ with $k\neq k'$: 
\begin{itemize}
\item If $k,k'$ aren't coprime, then $m(M_{i_{j+1},k}),m(M_{i_{j+1},k'})$ aren't coprime. 
\item For each $x,y \in M_{i_{j+1},k}$ there is $l \in \mathbb{N} \, \cup \{0\}$ and $a_1,\dots,a_l \in M_{i_{j+1},k}$ such that $\gcd(x,a_1),\gcd(a_1,a_2),\dots,\gcd(a_{l-1},a_l),\gcd(a_l,y)>1$.
\end{itemize}
\item We can partitionize  $M_{i_{j+1}}$ in $\{M'_{i_{j+1},k}\}_{k \in M_{i_j}}$ with the above conditions, only the second condition does not necessarily hold for one (arbitrary) $k_0 \in M_{i_j}$.
\end{itemize} 
\end{Corollary}

\section*{Further Remarks}
The fact, that we could use the Graph Minor Theorem here, has a simple background. Any structure, that could be reformulated in Graph Theory and with the possibility of a certain interpretation of minors, could be treated with the Graph Minor Theorem, at least infinite sequences of finite subsets of elements of it. Another interesting aspect is that there are some $S \in FCSG$ with the possibility that for every graph $G$ there is a subset $M$ of $S$, such that $G \simeq (M, \{ \{x,y\}\in [M]^2 | \gcd_S(x,y) \notin S^{\times} \})$, for example this is possible with $\mathbb{N}$. If it would be possible to prove the first part of Theorem \ref{the} for a certain $S \in FCSG$ with the mentioned trait without using the Graph Minor Theorem, we could get a shorter proof of it, because as remarked at the beginning of the paper the full proof was made in twenty papers and consisted of around 500 pages. Since nobody has reached to shorten it very much, it still remains important and interesting to get a short proof of one of the most important Theorems in Graph Theory.

\renewcommand{\refname}{References}

\vspace{0,3 cm}
\parindent=0pt
\end{document}